\documentclass[a4paper]{amsart}
\usepackage{amsmath}

\usepackage{amssymb}
\usepackage[all]{xy}
\usepackage{amsmath, amsfonts, amsthm , amssymb, amscd}
\usepackage{hyperref}
\usepackage{booktabs}
\usepackage{todonotes}
\usepackage{mathtools}
\usepackage{enumitem}
\usepackage{romannum}
\usepackage{anyfontsize}
\usepackage{tikz}
\usetikzlibrary{decorations.pathmorphing}
\newtheorem*{theorem*}{Theorem A}
\newtheorem*{theorem**}{Theorem B}

\theoremstyle{definition}

\theoremstyle{remark}

\title{The Magic Permutohedron}

\author{Djordje Barali\'{c} \and Lazar Milenkovi\'{c}}
\address{\scriptsize{ Mathematical Institute SASA, Knez Mihailova 36, p.p. 367, 11001, Belgrade, Serbia }}
\email{djbaralic@mi.sanu.ac.rs}
\address{\scriptsize{Tel Aviv University, P.O. Box 39040, Tel Aviv 6997801,  Tel Aviv, Israel}}
\email{milenkovic.lazar@gmail.com}

\subjclass[2020]{Primary 00A08. Secondary 05B30.}

\date{}

\begin{document}

\pagenumbering{arabic}

\begin{abstract}

We present a configuration called a magic permutohedron that shows the placement of the numbers of $\{1, 2, 3, \dots, 24\}$ in the vertices of the permutohedra so that the sum of numbers on each square side is 50 and the sum of the numbers in each hexagonal side is 75.

\end{abstract}

\maketitle

\section{What is the permutohedron?}

\begin{figure}[h!h!h!]
\centerline{\includegraphics[width=0.85\textwidth]{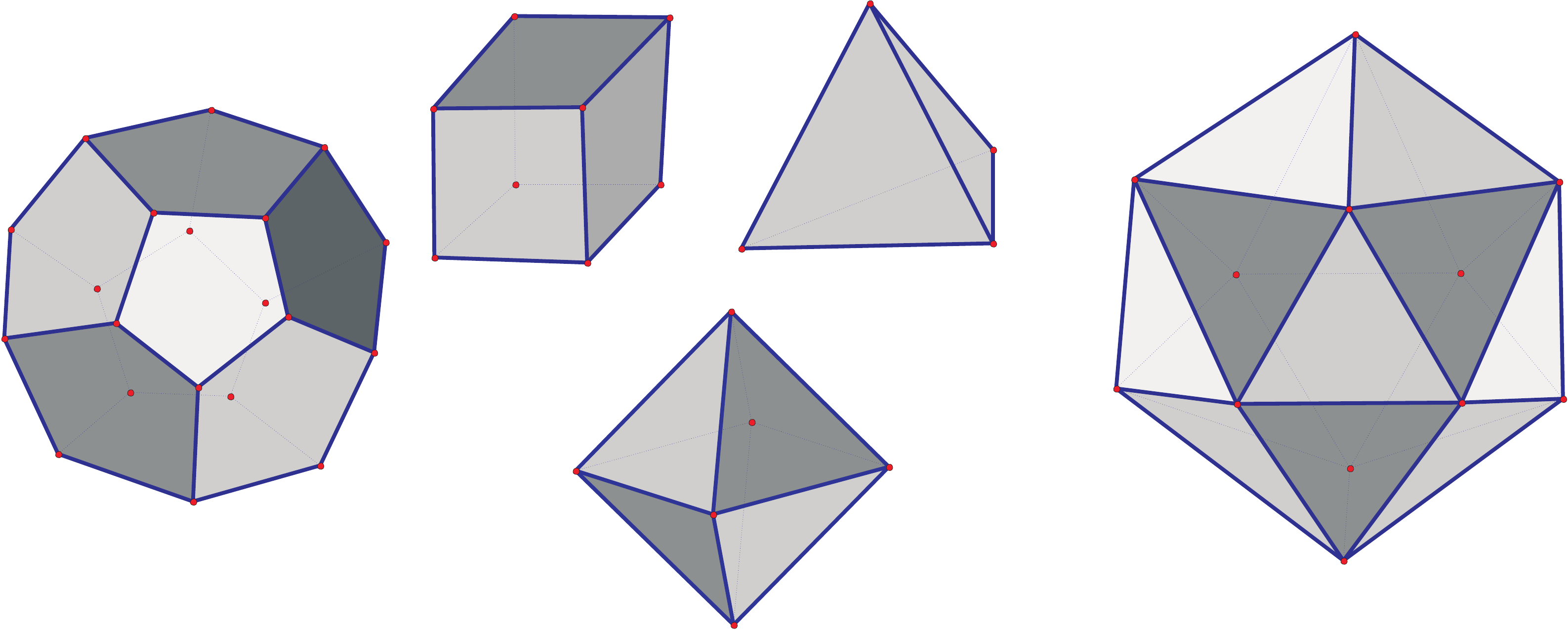}}
\caption{Platonic solids} \label{figurepla}
\end{figure}

Symmetric objects attracted human curiosity from ancient times. People often used to associate them with mystical and magical properties. However, they also motivated severe investigations and the birth of new mathematical and science disciplines and multidisciplinary studies. Early examples are five regular solids, also known as Platonic solids, see Figure~\ref{figurepla}.

Why are such objects so specific and appear in many fascinating natural places? Philosophically speaking, because they are perfect geometric objects and natural shapes tends to be perfect. Mathematics is a language we are trying to explain this perfectness. If geometry catches their forms and obvious exterior beauty, other mathematical disciplines like group theory, number theory and algebra describe their hidden deep properties.

\begin{figure}[h!h!h!]
\centerline{\includegraphics[width=\textwidth]{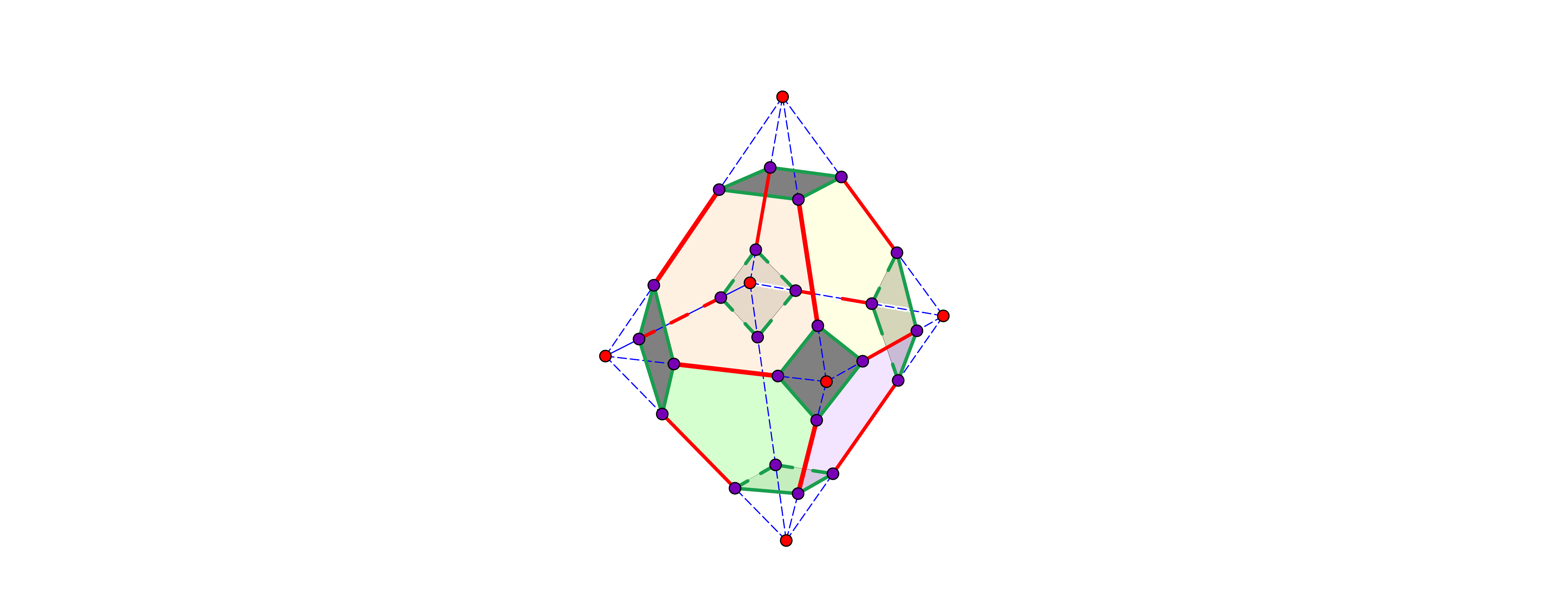}}
\caption{Truncated octahedra} \label{Figuretru}
\end{figure}

The solid we will present here is one of the essential objects in modern mathematics. It is obtained by truncating the regular octahedron in all its vertices, see Figure~\ref{Figuretru}. Its facets are regular polygons: six squares and eight regular hexagons. The solid has 24 vertices and 36 edges, with three edges meeting in each vertex. Its first recorded appearance \cite{Schoute} in mathematical literature dates back to 1911 and a Dutch mathematician Pieter Hendrik Schoute known for his work on regular polyhedrons and their higher dimension analogues, \textit{regular polytopes}. Georges Th. Guilbaud and Pierre Rosenstiehl~\cite{Rosen} named it \textit{the permutohedron} in 1963.

Its name is coined from two words, permute and hedra (the face of a solid). There are 24 permutations of the elements of the set $\{1, 2, 3, 4\}$ and 24 vertices of the permutohedron, but could we give a more satisfactory explanation to justify its name? Of course, the answer is yes! The simplest way is to introduce it as a member of a family of convex polytopes. A convex polytope in a $n$-dimensional Euclidean space $\mathbb{R}^n$ is a convex hull of a finite set of points in this space. We maintained the basic features of convex polygons and polyhedra for higher dimensions by this definition. Indeed, geometrical analogues between polygons and polyhedra in dimensions two and three emerge here. For example, a convex $4$ polytope beside vertices and edges will have polygons as $2$-faces  and polyhedra as its $3$-faces. More generally, a convex $n$-polytope has faces of dimension $k$ for all $0\leq k\leq n-1$, while faces of dimension $n-1$ are called facets. The interested reader in polytopes is kindly referred to the outstanding book \textit{Lectures on Polytopes} by G\"{u}nter Ziegler \cite{Zieg}.

\begin{figure}[h!h!h!]
\centerline{\includegraphics[width=0.85\textwidth]{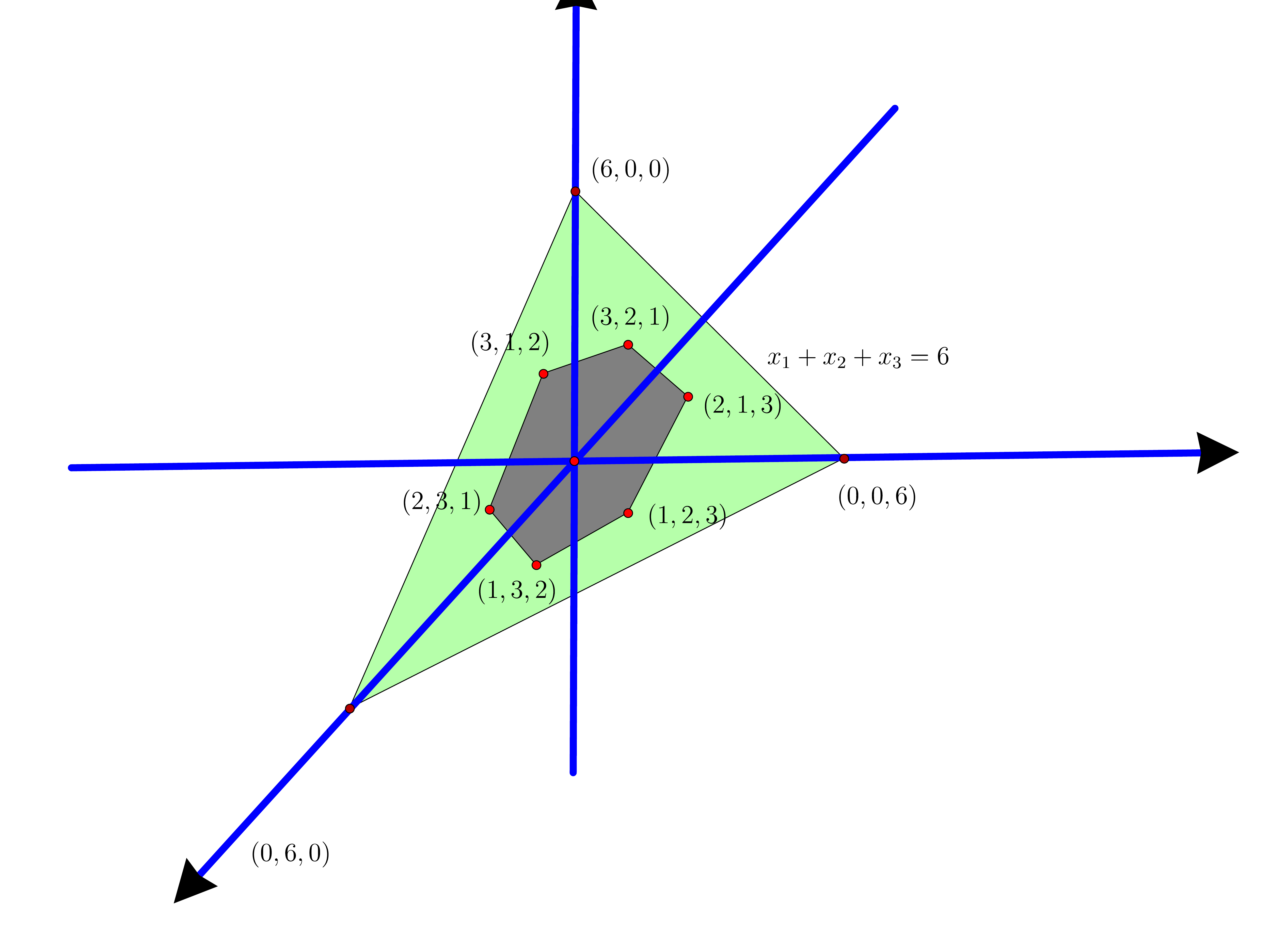}}
\caption{Permutohedron as the convex hull of its vertices} \label{figurehex}
\end{figure}

\begin{figure}[h!h!h!]
\centerline{\includegraphics[width=\textwidth]{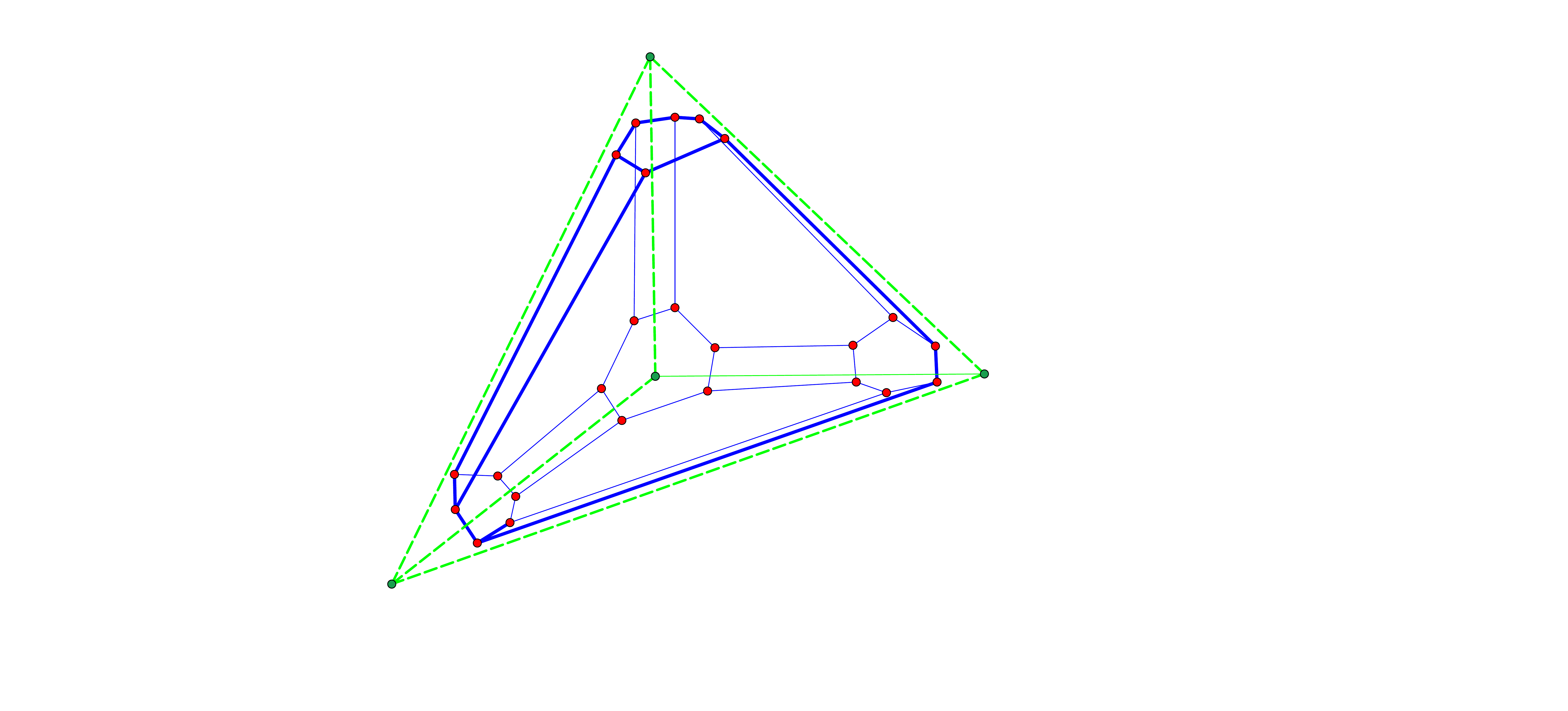}}
\caption{Total truncation of tetrahedra} \label{figuretrun}
\end{figure}

The $n$-permutohedron is the convex hull of $(n+1)!$ points in $\mathbb{R}^{n+1}$, each one of them being a permutation of the coordinates of $\{1, 2, 3, \dots, n\}$. A regular hexagon is therefore a 2-dimensional permutohedron. As the sum of vertex coordinates is always $\frac{(n+1)(n+2)}{2}$, the $n$-permutohedron lies in the hyperplane $$x_1+x_2+\cdots+x_n=\frac{(n+1)(n+2)}{2}$$ in $\mathbb{R}^{n+1}$, see Figure~\ref{figurehex}. It has between $\frac{n(n+1)!}{2}$ and $2^{n+1}-2$ facets, while each facet corresponds to a proper subset of $\{1, 2, \dots, n+1\}$. The last fact can be seen geometrically because the $n$-permutohedron can be defined alternatively after truncating all faces of $n$ simplex, see Figure~\ref{figuretrun}. An $n$-simplex is a convex hull of $n+1$ points in general positions in $\mathbb{R}^n$, a $n$-dimensional analogue of a triangle and a tetrahedra. The combinatorics around $k$-faces of $n$-permutohedron is exciting. The total number of faces of dimension $k$ is equal to $(n-k+1) \cdot S (n+1, n+1-k)$, where $S (n, k)$ stands for the Stirling number of the second kind - the number of ways to partition a set of $n$ elements into $k$ non-empty subsets.

\begin{figure}[h!h!h!]
\centerline{\includegraphics[width=0.85\textwidth]{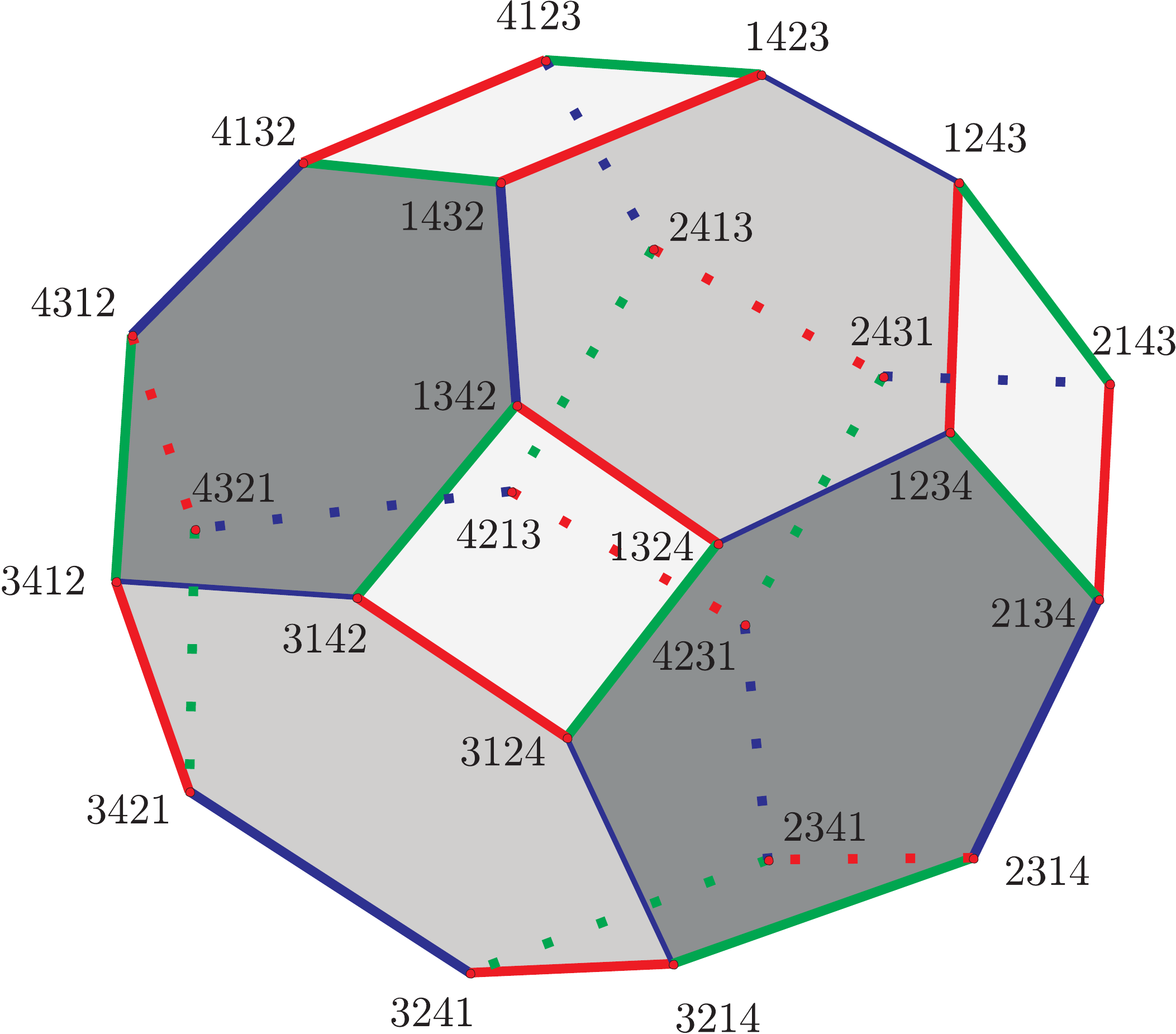}}
\caption{Permutohedra as a Cayley graph of $S_4$} \label{fig}
\end{figure}

The permutation group $S_{n+1}$ acts on the vertices of the $n$-permutohedron by permuting the coordinates, and this is the principal reason this object has a fancy look. Indeed, the permutohedron helps understand permutations and the symmetric group $S_{n+1}$. To any group $G$ and a set of its generators $S$, we can assign an edge-colored directed graph $\Gamma (G, S)$ whose vertex set are elements of $G$ and there is a directed edge colored by $s\in S$ from $g\in G$ to $h\in G$ exactly when $h=g s$ in $G$. Such a graph is called a Cayley graph of $G$. It turns out that the vertices of the $n$-permutohedron can be labeled by permutations so that its edges make the Cayley graph of $S_{n+1}$ for the set of generators being transpositions $$(1,2), (2, 3), \dots, (n, n+1), $$ those that swap two consecutive elements, see Figure~\ref{fig}.

\section{Magic figures and solids}

\begin{figure}[h!h!h!]
\centerline{\includegraphics[width=\textwidth]{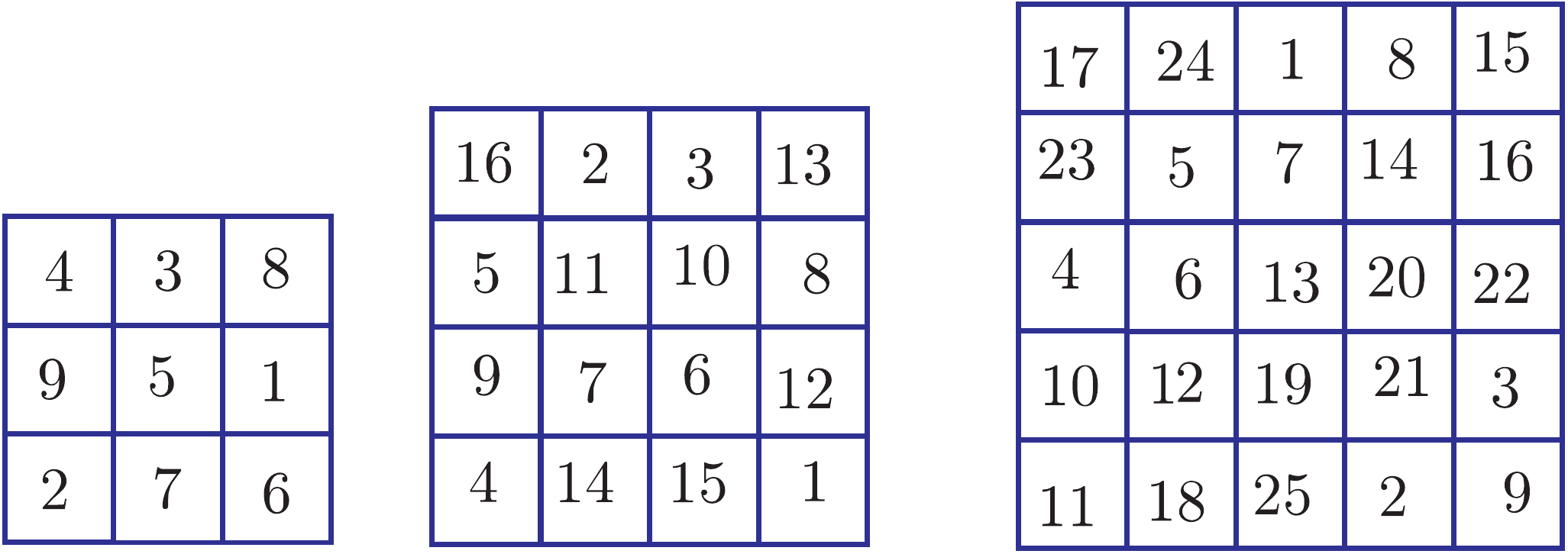}}
\caption{Magic squares} \label{magic}
\end{figure}

The magic square is one of the earliest recreational mathematical objects known at least 190 BCE. People assigned them an occult or mythical significance in the past, so they often appeared as symbols in works of art. The magic square of order $n$ is a sequence of $n^2$ numbers arranged in $n\times n$ squares so that the sum of numbers in each row, each column and both diagonals are the same, see Figure~\ref{magic}. The latter number is also known as a \textit{magic constant}. Mathematicians also studied magic squares with extra constraints, such as ultra magic squares, pandiagonal, symmetric magic squares, etc.

Depending on the remainder of $n$ modulo 4 there are various methods for constructing higher-order magic squares. If $n$ is odd, we say that a magic square is odd. It is called doubly even if four divides $n$ and singly even when $n \equiv 2 \pmod 4$. Singly even magic squares are more challenging to be generated. Some methods allow us to create more than one magic square for a given $n$, while others, which are usually simpler, provide us with just one magic square of a given order. Some transformations and symmetries allow us to construct a few more solutions. However, classifying all magic squares of a given order has been an open question for a long time. The problem has been studied by many famous mathematicians, including John Conway, Simon de la Loub\`{e}re, and Claude Bachet, who wrote  \textit{Probl\`{e}mes plaisans et delectables qui se font par les nombres},  the first books in the recreational mathematics in early 17th century.

\begin{figure}[h!h!h!]
\centerline{\includegraphics[width=0.45\textwidth]{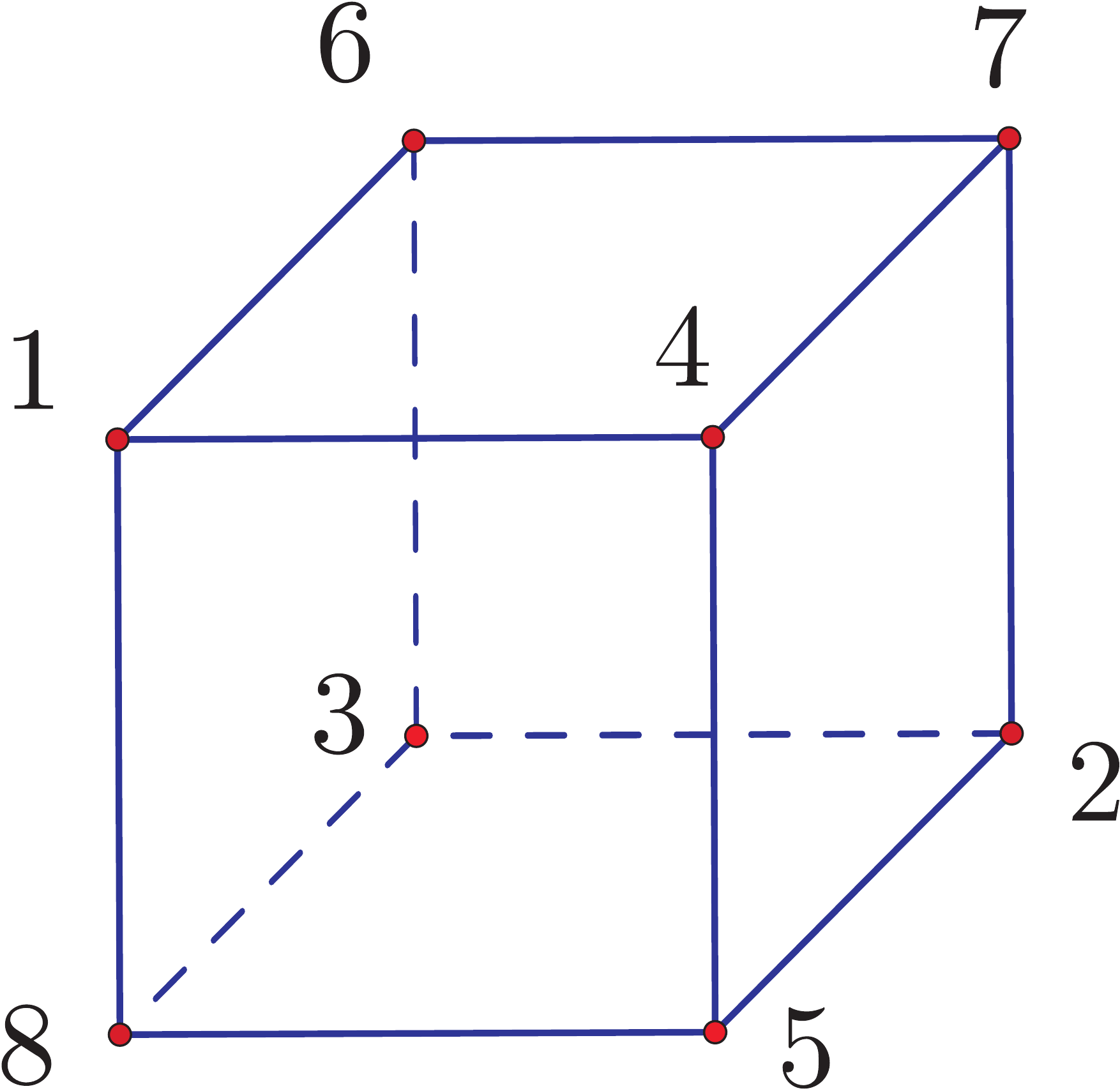}}
\caption{Magic cube} \label{magiccube}
\end{figure}

The idea of magic squares was implemented in other plane figures and even in higher dimensional figures. Here we consider a similar problem where the goal is arranging the numbers from the set of $\{1, 2, \dots, m\}$ in the vertices of a regular $m$-vertex solid. We immediately ruled out such possibilities for the tetrahedron and the icosahedron since two triangles sharing an edge would have had the same numbers in the vertices opposite to the joint edge. Finding such a placement is straightforward; see \ref {magiccube} for an example. It is easy to see that in any such configuration, the sum of the vertices on each face of the cube has to be 18. A small analysis shows six distinct options up to the cube's rotations.

Let us think about the dodecahedron case. If we denote by $S$ the sum of numbers in each pentagon, then we have that $$12 S=3 \cdot (1+2+\cdots+20)=630,$$ but this is impossible since $12\nmid 630$. Apart from the cube, other Platonic solids do not have surprising configurations for placing the numbers in the vertices. One can ask what is happening with other polyhedra and whether we can relax the constraint so that we still have some magical properties.

\begin{figure}[h!h!h!]
\centerline{\includegraphics[width=0.45\textwidth]{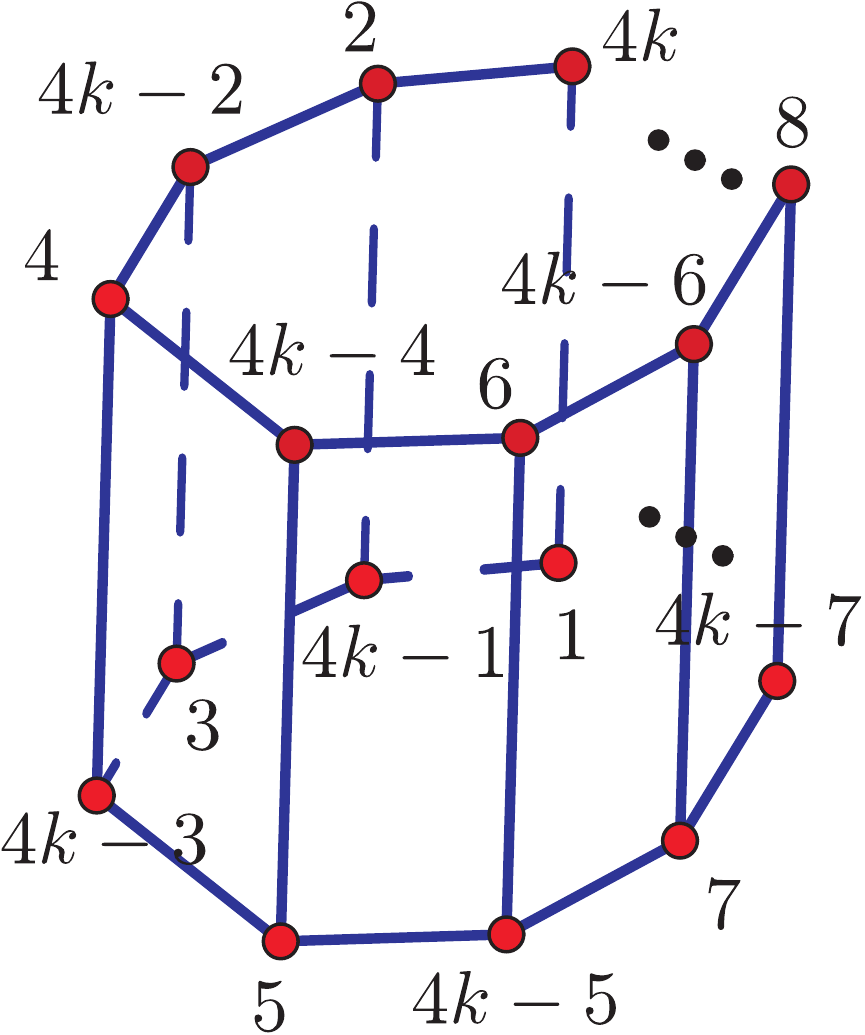}}
\caption{Magic $2k$-gon prism} \label{prism}
\end{figure}

The next example of a regular $n$-gon prism is illustrative. Let us try to arrange $1$, $2, \dots, 2n$ in the vertices of this prism so that the sum of the numbers in each basis is $A$, and the sum of the numbers in each rectangular face of the prism is $B$. We have the following conditions:
\begin{align*}
2 A= \frac{2n (2n+1)}{2}=n(2n+1) \label{eq1} & \text{\, \, and \, \,} n B= 2n (2 n+1)\\
\end{align*}
The second condition yields $B=4n+2$, but the requirement from the first condition is that $n$ must be even. It is easy to produce such a configuration for $n=2k$; one possibility is depicted in Figure~\ref {prism}.

\section{Magic permutohedron}

The example of the magic $2n$-gon prism inspires us to consider the following question:

\textit{Is it possible to place the numbers $1, 2, \dots, 24$ in the vertices of permutohedra so that the sum of the numbers in each square face is the same and the sum of the numbers in each hexagonal face is also the same?}

We straightforwardly find that the sum of the numbers in a square face and a hexagonal face should be 50 and 75, respectively. However, if we try to find an arrangement by hand explicitly, we are most likely to fail as we have too many spaces of freedom at the beginning. On the other hand, the authors' attempts to put some more natural constraints failed, so we wrote a program to make it for us!

\begin{figure}[h!h!h!]
\centerline{\includegraphics[width=0.65\textwidth]{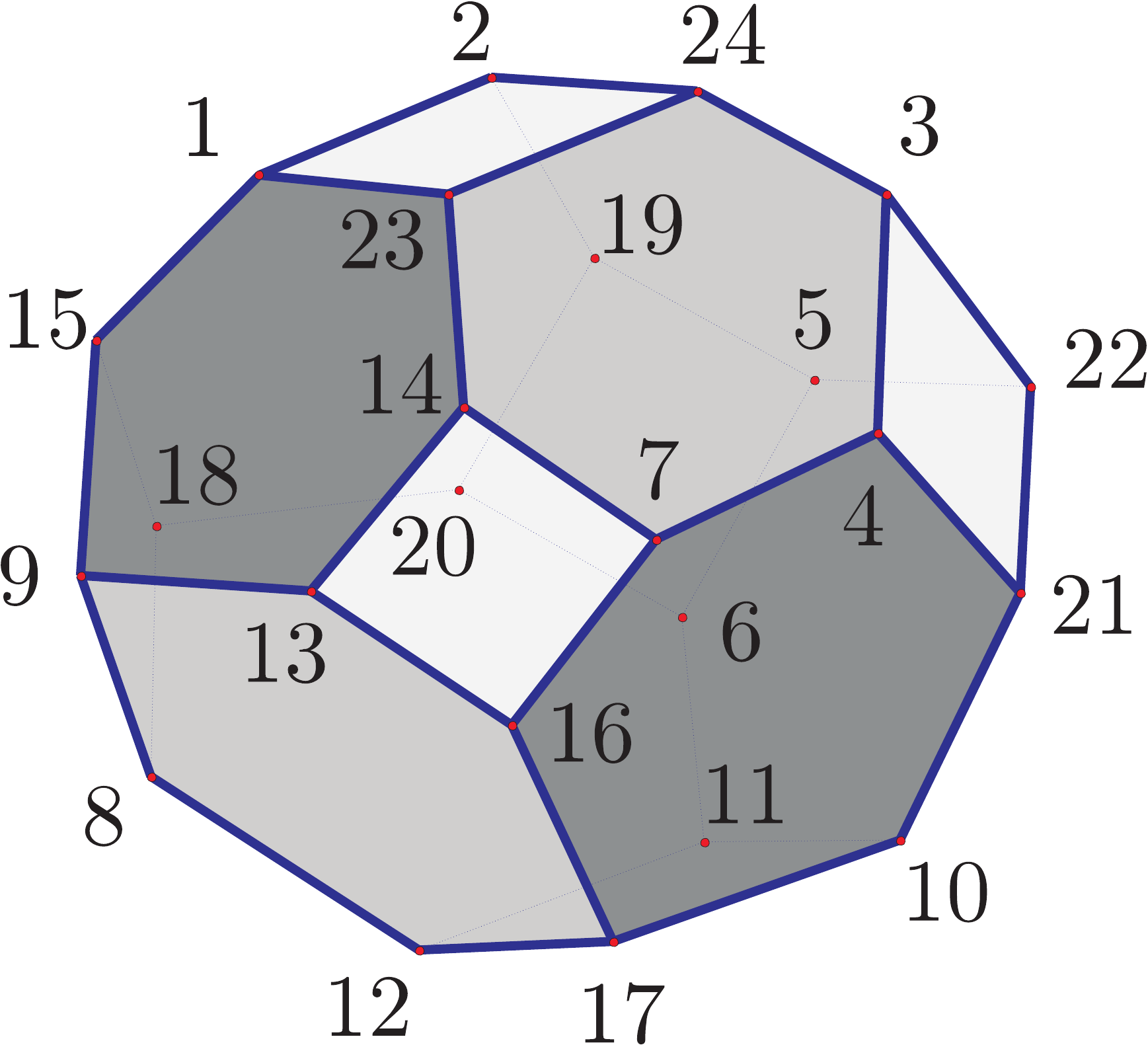}}
\caption{The magic permutohedra} \label{mp}
\end{figure}

The results were more than exciting. The computer found $3900064$ distinct solutions (up to the rotations of the permutohedra). However, this number of solutions makes less than $4\time 10^{-15} \%$ of the total of 12926008369442488320000 possibilities. One solution found by the computer is presented in Figure~\ref{mp}.

We wonder if a similar property holds for the $n$-permutohedron for $n\geq 4$. The condition can be formulated on 2-faces (hexagons and squares) or on its facets which are the products of smaller dimensional permutohedra.

\section*{Acknowledgements}

The first author was supported by the Serbian Ministry of Science, Innovations and Technological Development through the Mathematical Institute of the Serbian Academy of Sciences and Arts.

\end{document}